\documentclass[12pt,a4paper,reqno]{amsart}
\usepackage{amsmath,amscd,amssymb,latexsym}
\usepackage{longtable}

\newtheorem{Thm}{Theorem}[section]

\newtheorem{Lem}[Thm]{Lemma}
\newtheorem{Prop}[Thm]{Proposition}
\theoremstyle{definition}

\newtheorem{Exm}[Thm]{Example}
\theoremstyle{remark}
\newtheorem{Rk}[Thm]{Remark}

\newcommand{\Z}{\mathbb{Z}}

\newcommand{\PP}{\mathbb{P}}
\newcommand{\md}{\,\operatorname{mod}}
\newcommand{\sm}{\left(\smallmatrix}
\newcommand{\esm}{\endsmallmatrix\right)}

\begin{document}
\title{On the arithmetic of certain modular curves}
\author[D. Jeon \and C. H. Kim]{Daeyeol Jeon \and Chang Heon Kim}

\address{Daeyeol Jeon, Department of Mathematics Education, Kongju National University,
182 Shinkwan-dong, Kongju, Chungnam, 314-701 Korea}
\email{dyjeon@kongju.ac.kr}

\address{Chang Heon Kim, Department of mathematics, Seoul Women's university,
126 Kongnung 2-dong, Nowon-gu, Seoul, 139-774 Korea}
\email{chkim@swu.ac.kr}


\begin{abstract}
In this work, we estimate the genus of the intermediate modular
curves between $X_1(N)$ and $X_0(N),$ and determine all the
trigonal ones.
\end{abstract}

\maketitle \setcounter{section}{-1}

\section{Introduction}
Let $N$ be a positive integer and $\Delta$ a subgroup of
$({\mathbb Z}/{N \mathbb Z})^*$ which contains $\pm1.$ Let
$X_\Delta(N)$ be the modular curve defined over $\mathbb Q$
associated to the congruence subgroup
$$\Gamma_\Delta(N):=\left\{{\begin{pmatrix}a&b\\c&d\end{pmatrix}}
\in\operatorname{SL}_2(\Z)\,|\,\,a\md N\in\Delta, N\mid c
\right\}.$$ Then all the intermediate modular curves between
$X_1(N)$ and $X_0(N)$ are of the form $X_\Delta(N).$ Denote the
genus of $X_\Delta(N)$ by $g_\Delta(N).$ In this paper we carry
out some arithmetic of the curves $X_\Delta(N).$

In Section $1$ we prove a genus formula of the curves
$X_\Delta(N)$ which was referred in the authors' previous works
\cite{J-K,J-K2,J-K-S} without proof.

A smooth projective curve $X$ defined over an algebraically closed
field $k$ is called {\it $d$-gonal} if it admits a map $\phi:X\to
\PP^1$ over $k$ of degree $d$. For $d=3$ we say that the curve is
$trigonal.$ Also, the smallest possible $d$ is called the {\it
gonality} of the curve and denote it by $\operatorname{Gon}(X).$

Hasegawa and Shimura \cite{H-S} proved that $X_0(N)$ is trigonal
if and only if it is of genus $g\leq 2$ or is not hyperelliptic of
genus $g=3,4.$ In fact the ``if"-part is well-known. The modular
curves $X_0(N)$ carry the action of the Atkin-Lehner involutions
$W_d$ for any $d\|N$ which denote a positive integer $d$ dividing
$N$ with $(d,N/d)=1$. Let $X_0^{+d}(N)$ and $X_0^{*}(N)$ be the
quotient of $X_0(N)$ by $W_d$ and the $W_d$'s for all $d\|N$
respectively. In \cite{H-S2,H-S3}, they determined also the
trigonal modular curves $X_0^{+d}(N)$ and $X_0^{*}(N),$ and found
that there exist non-trivial trigonal modular curves, i.e., those
of genus $g\geq 5.$

The authors and Schweizer \cite{J-K-S} showed that there exist no
non-trivial trigonal modular curves $X_1(N),$ which plays a
central role in the determining the torsion structures of elliptic
curves defined over cubic number fields which occur infinitely
often.

In Section $3$ we determine all the intermediate modular curves
between $X_1(N)$ and $X_0(N)$ which are trigonal, and conclude
that there exist no non-trivial trigonal curves. For the purpose,
it is necessary to determine all the hyperelliptic intermediate
modular curves, which was done by Ishii and Momose \cite{I-M}. But
there was a mistake, and so we complete their results in Section
$2.$

\section{A genus formula}
Let $\Gamma(1)=\operatorname{SL}_2(\mathbb Z)$ be the full modular
group. For any integer $N\geq 1,$ we have subgroups $\Gamma_1(N)$
and $\Gamma_0(N)$ of $\Gamma(1)$ defined by matrices $\sm
a&b\\c&d\esm$ congruent modulo $N$ to $\sm 1&*\\0&1\esm$ and $\sm
*&*\\0&*\esm$ respectively. We let $X_1(N)$ and $X_0(N)$ be the
modular curves defined over $\mathbb Q$ associated to
$\Gamma_1(N)$ and $\Gamma_0(N)$ respectively. The $X$'s are
compact Riemann surfaces.  Let $g_0(N)$ denote the genus of
$X_0(N).$ For any congruence subgroup $\Gamma\subset\Gamma(1),$ we
shall denote by $\overline\Gamma$ the image of $\Gamma$ by the
natural map $\Gamma(1)\to\overline\Gamma(1):=\Gamma(1)/\{\pm1\}.$

For $d|N,$ let $\pi_d$ be the natural projection from $({\mathbb
Z}/{N \mathbb Z})^*$ to $({\mathbb Z}/{\{d,N/d\} \mathbb Z})^*$,
where $\{d,N/d\}$ is the least common multiple of $d$ and $N/d.$
Then we have the following genus formula:

\begin{Thm}\label{genus1}
The genus of the modular curve $X_\Delta(N)$ is given by
$$g_\Delta(N)=1+\frac{\mu}{12}-\frac{\nu_2}{4}-\frac{\nu_3}3-\frac{\nu_\infty}2$$
where
\begin{eqnarray*}
\mu &=& N\prod_{p|N\atop
prime}\left(1+\frac1p\right)\frac{\varphi(N)}{|\Delta|}\\
\nu_2 &=& \left|\{(b\md N)\in\Delta\,|\, b^2+1\equiv 0\md
N\}\right|\frac{\varphi(N)}{|\Delta|}\\ \nu_3 &=& \left|\{(b\md
N)\in\Delta\,|\, b^2-b+1\equiv 0\md
N\}\right|\frac{\varphi(N)}{|\Delta|}\\ \nu_\infty &=&
\sum_{d|N\atop d>0}\frac{\varphi(d)\varphi(\frac Nd
)}{|\pi_d(\Delta)|}.\end{eqnarray*}
\end{Thm}
\begin{proof}
We follow the notations of \cite{Ogg1}. One has to check that the
index of $\overline{\Gamma}_\Delta(N)$ in $\overline{\Gamma}(1)$
is $\mu$, the number of elliptic fixed points of order $2$ (resp.
$3$) is $\nu_2$ (resp. $\nu_3$) and that the number of cusps is
$\nu_\infty$. It is easy to show the following:
\begin{align*}
\mu & =[\overline{\Gamma}(1): \overline{\Gamma}_\Delta(N)]
=[\overline{\Gamma}(1):\overline{\Gamma}_0(N)]
[\overline{\Gamma}_0(N): \overline{\Gamma}_\Delta(N)] \\
    & =N \prod_{p|N}\left(1+\frac1p\right)\frac{\varphi(N)}{|\Delta|}.
\end{align*}

Put $L_0=\sm N&0\\0&1 \esm$. Then the double coset $\Gamma(1) L_0
\Gamma(1)$ has the right coset decomposition as follows:
$$\Gamma(1) L_0 \Gamma(1)=\bigcup \Gamma(1) L$$
where $L=\sm a&b\\0&d \esm$ with $a>0$, $ad=N$, $b \md d$ and
$(a,b,d)=1$.

Now we compute $\nu_2$ and $\nu_3.$ Let $A$ be an elliptic element
in $\Gamma(1)$ and $P$ the fixed point of $A$ in the complex upper
half plane. Then $P=M P_0$ for some $M\in \Gamma(1)$ where $P_0=i$
or $e^{{2\pi i}/3}$. Write $L_0 M= B L$ for some $B\in \Gamma(1)$
and $L=\sm a&b\\0&d \esm$ with $a>0$, $ad=N$ and $(a,b,d)=1$. Now
if $P_0=i$,
\begin{align*}
 & A=M\sm 0&-1\\1&0 \esm M^{-1} \in \Gamma_0(N)
 =\Gamma(1) \cap L_0^{-1} \Gamma(1) L_0\\
 & \iff L \sm 0&-1\\1&0 \esm L^{-1} \in \Gamma(1) \\
 & \iff \sm a&b\\0&d \esm \sm 0&-1\\1&0 \esm \sm a&b\\0&d
 \esm^{-1}= \sm b/a & -(a^2+b^2)/N \\ d/a & b/a \esm\in \Gamma(1)
 \\
 & \iff a=1, d=N \text{ and } b^2+1 \equiv 0 \md N.
\end{align*}
Similarly if $P_0=e^{2\pi i}/3$,
 \begin{align*}
 A\in \Gamma_0(N)
 & \iff \sm a&b\\0&d \esm \sm 0&-1\\1&1 \esm \sm a&b\\0&d
 \esm^{-1}= \sm b/a & -(a^2-ab+b^2)/N \\ d/a & (a-b)/a \esm\in \Gamma(1) \\
 & \iff a=1, d=N \text{ and } b^2-b+1 \equiv 0 \md N.
\end{align*}
 Write $M=\sm x&y \\ z&w \esm$ and $B=\sm x'&y' \\ z'&w' \esm$.
 From $L_0 M=BL$ it follows that
 \begin{equation} \label{AA}
 \sm  Nx&Ny \\ z&w \esm = \sm
  x'& bx'+Ny'\\ z'& bz'+Nw' \esm.
 \end{equation}
Note that $M \sm 0&-1\\1&0 \esm M^{-1}=\sm yw+xz&*\\ *&* \esm$ and
$M \sm 0&-1\\1&1 \esm M^{-1}=\sm yw+xz-yz&*\\ *&* \esm.$ Then for
the elliptic element $A$ of order $2$ (resp. $3$) to lie in
$\Gamma_\Delta(N)$, we need $yw+xz \md N \in \Delta$ (resp.
$yw+xz-yz \md N \in \Delta$) together with the condition
$b^2+1\equiv 0 \md N$ (resp. $b^2-b+1\equiv 0 \md N$). From
(\ref{AA}) it is easy to see that $yw+xz \equiv -b \md N$ and
$yw+xz-yz \equiv -b+1 \md N$. Thus if $A$ is an elliptic element
of order $2$ (resp. $3$) in $\overline{\Gamma}_\Delta$, then there
corresponds an element $b \md N \in \Delta$ satisfying
$b^2+1\equiv 0 \md N$ (resp. $b^2-b+1\equiv 0 \md N$). Conversely,
we can form an elliptic element of order 2 (resp. 3) from the
solution in $\Delta$ of the congruence equation $x^2+1 \equiv 0
\md N$ (resp. $x^2-x+1 \equiv 0 \md N$). We note that different
solutions give $\Gamma_0(N)$-inequivalent elliptic points of order
$2$ (resp. $3$).

Now we consider the Galois covering $p_2: X_\Delta(N) \to X_0(N)$.
If $A$ is an elliptic element of order $2$ in
$\overline{\Gamma}_\Delta$ and $AP=P$, then each point in the
inverse image of $\Gamma_0(N)P$ is again an elliptic point of
order $2$ and has ramification index $1$. Thus the number of
elements in $p_2^{-1} (\Gamma_0(N)P)$ would become the degree of
$p_2,$ and hence we have the following:
 \begin{align*}
 \nu_2 & =
 \left|\{(b\md N)\in\Delta\,|\, b^2+1\equiv 0\md
  N\}\right|\cdot \text{degree of } p_2 \\
  & = \left|\{(b\md N)\in\Delta\,|\, b^2+1\equiv 0\md
    N\}\right|\frac{\varphi(N)}{|\Delta|}.
 \end{align*}
 Similarly, $\nu_3=\left|\{(b\md N)\in\Delta\,|\, b^2-b+1\equiv 0\md
    N\}\right|\frac{\varphi(N)}{|\Delta|}.$

Finally we compute $\nu_\infty.$ We follow the notations of
\cite{Ogg2}. Let $p_1: X_1(N) \to X_\Delta(N)$, $p_2: X_\Delta(N)
\to X_0(N)$ be the Galois coverings and $p=p_2 \circ p_1$. Denote by $s=\sm x \\
y \esm$ a cusp in $X_1(N)$. Then the following holds:
$$e_p(s)=e_{p_1}(s) e_{p_2}(p_1 s) \text{ and } e_p(s)=(N/d,d)
\text{ with } d=(y,N)$$ where $e$'s denote ramification indices
(\cite{Ogg2}, Proposition 2). Now we claim that
$$ e_{p_1}(s)=|\Delta|/|\pi_d (\Delta)|. $$

Note that the group $G=\Gamma_\Delta(N)/\pm\Gamma_1(N)$ is
isomorphic to $\Delta/\pm 1.$  Each element $\sm a&b\\c&d \esm\in
\Gamma_\Delta(N)$ acts on $\sm x\\y\esm$ as $\sm ax\\a^{-1}y\esm.$
Then $\sm x\\y\esm$ and $\sm ax\\a^{-1}y\esm$ represent the same
cusp on $X_1(N)$ if and only if $ax\equiv \pm x\md d$ and
$ay\equiv \pm y\md N$, i.e., $a\equiv \pm 1\md d$ and $\frac{N}d.$

Let $\{d,\frac{N}d\}$ denote the least common multiple of $d$ and
$\frac{N}{d}.$ Let $H=\bigl\{a\md N \in \Delta/ \pm 1 \, | \,
a\equiv 1\md\{d,\frac{N}d\}\bigl\}.$ Since $H$ is the kernel of
the natural map $\Delta/ \pm 1 \to (\mathbb
Z/\{d,\frac{N}d\})^*/\pm 1,$ the number of $H$ is equal to
$|\Delta|/|\pi_d(\Delta)|.$ We can view $H$ as a subgroup of $G.$
Then $G/H$ has the same cardinality as the set of orbits $Gs.$
Since the elements of $Gs$ are the cusps in $X_1(N)$ lying over
the cusp $p_1(s)$ in $X_\Delta(N),$ the ramification index of $s$
in $X_1(N)$ is equal to the cardinality of $H$. By the claim we
come up with
 \begin{align*}
 \nu_\infty
 & = \sum_{d|N \atop d>0} \frac{\deg p_2}{e_{p_2}}
     \varphi((d,N/d)) \text{ \, since $p_2$ is a Galois covering }
     \\
 & = \sum_{d|N \atop d>0} \frac{\varphi (N)}{|\Delta|}
     \frac{|\Delta|}{|\pi_d(\Delta)|} \frac{1}{(N/d,d)}
     \varphi((d,N/d))
     \\
 & = \sum_{d|N \atop d>0} \varphi(d) \varphi(N/d)/|\pi_d(\Delta)|.
\end{align*}
The last equality can be shown by using the fact that
$$\varphi(n_1)\varphi(n_2)=\varphi(n_1 n_2)
\frac{\varphi((n_1,n_2))}{(n_1,n_2)}.$$
\end{proof}

\section{hyperelliptic modular curves}
If a curve $X$ is $2$-gonal, we call it {\it sub-hyperelliptic.}
Also if $X$ is sub-hyperelliptic of genus $g\geq 2,$ then it is
called {\it hyperelliptic}.

\begin{Prop}\cite{Ne,N-S}\label{d-gonal}
Let $X_1$ and $X_2$ be smooth projective curves over an
algebraically closed field $k,$ and assume that there is a finite
morphism $X_1\to X_2$ over $k.$ If $X_1$ is $d$-gonal, so is
$X_2.$
\end{Prop}

The best general lower bound for the gonality of a modular curve
seems to be the one that is obtained in the following way.

Let $\lambda_1$ be the smallest positive eigenvalue of the
Laplacian operator on the Hilbert space $L^2 (X_\Gamma)$ where
$X_\Gamma$ is the modular curve corresponding to a congruence
subgroup $\Gamma$ of $\Gamma(1)$ and $D_\Gamma$ be the index of
$\overline\Gamma$ in $\overline\Gamma(1).$ Abramovich \cite{A}
shows the following inequality:
$$\lambda_1 D_\Gamma \le 24 \operatorname{Gon}(X_\Gamma).$$
Using the best known lower bound for $\lambda_1$, due to Henry
Kim and Peter Sarnak, as reported on page 18 of \cite{B-G-G-P},
i.e., $\lambda_1> 0.238$, we get the following result.

\begin{Thm}\label{bound1}
Let $X_\Gamma$ be the modular curve corresponding to a congruence
subgroup $\Gamma$ of index
$D_\Gamma:=[\overline\Gamma(1):\overline\Gamma].$ Then
$$D_\Gamma < \frac{12000}{119}\operatorname{Gon}(X_\Gamma).$$
\end{Thm}

In the following, we call the inequality in Theorem \ref{bound1}
Abramovich's bound.

Ishii and Momose \cite{I-M} asserted that there exist no
hyperelliptic modular curves $X_\Delta(N)$ with $\{\pm
1\}\varsubsetneq\Delta\varsubsetneq(\Z/N\Z)^*$. But we get the
following result.

\begin{Thm}\label{hyperelliptic}
There exists one and only  hyperelliptic modular curve of the form
$X_\Delta(N)$ with $\{\pm
1\}\varsubsetneq\Delta\varsubsetneq(\Z/N\Z)^*$, i.e.,
$X_{\Delta_1}(21)$ where $\Delta_1$ is in Table \ref{lowgenus}.
\end{Thm}

\begin{Rk}
In \cite{I-M} there was a mistake to treat Atkin-Lehner
involutions on $X_\Delta(N).$ The Atkin-Lehner involutions define
a unique involution on $X_0(N)$ but it doesn't hold for
$X_\Delta(N).$
\end{Rk}

From now on we prove Theorem \ref{hyperelliptic}. To this end we
need some preparations.

Let $X$ be a smooth projective curve of genus $g\geq 2$ and
$\Omega^1(X)$ the space of holomorphic differential forms on $X.$
Then $\Omega^1(X)$ gives rise to a line bundle, and this is called
the {\it canonical bundle}. This line bundle, in certain
situations, gives an embedding into projective space.

Let $\omega_1\dots,\omega_g$ be a basis for $\Omega^1(X)$. Viewing
$X$ as a Riemann surface, we may choose a finite covering of open
sets, with local parameters $z$ on each set, such that we can
locally write $\omega_i=f_i(z)dz.$ Then we get the following
well-defined map
$$\begin{array}{llll}
\phi: & X & \to & \PP^{g-1}\\
& P & \mapsto & (\omega_1(P):\dots :\omega_g(P)).
\end{array}$$
Note that $(\omega_1(P):\dots :\omega_g(P))=(f_1(P):\dots
:f_g(P)).$ The above map is called a {\it canonical map}. Let
$\overline X$ denote the image of $X$ by the canonical map. It is
well-known that if $X$ is not hyperelliptic then the canonical map
is injective.

If one takes the canonical map of a hyperelliptic curve then it is
possible to predict what will happen. Suppose $X$ is a
hyperelliptic curve of genus $g\geq 3.$ Then the image $\overline
X$ under the canonical map will be a smooth curve which is
isomorphic to $\PP^1$ and which is described by $(g-1)(g-2)/2$
quadratic equations (See \S 2 of \cite{Ga}).

Therefore it is possible to distinguish between hyperelliptic
curves and non-hyperelliptic ones by examining their images under
the canonical map.

Now we consider the modular curves $X_\Delta(N)$ of genus
$g=g_\Delta(N)\geq 3.$ Let $S_\Delta^2(N)$ denote the space of
cusp forms of weight $2.$ Suppose $\{f_1,\dots,f_g\}$ is a basis
of $S_\Delta^2(N)$. Then the canonical map may be written as
$$X_\Delta(N)\ni P\mapsto (f_1(P):\cdots:f_g(P))\in{\mathbb
P}^{g-1}.$$ One can get such a basis and their Fourier
coefficients from \cite{St}. Then to obtain a system of quadratic
generators of $I(\overline{X_\Delta(N)}),$ we have only to compute
the relations of the $f_if_j\,(1\leq i,j\leq g).$ If $X_\Delta(N)$
is not hyperelliptic, then there exist exactly $(g-2)(g-3)/2$
linear relations among the $f_if_j$ (See \S 2 of \cite{H-S}).

Now we are ready to prove Theorem \ref{hyperelliptic}. By
Proposition \ref{d-gonal} it suffices to consider $X_\Delta(N)$
when $X_0(N)$ are sub-hyperelliptic. If the genus $g_0(N)\leq 2,$
then one can find all $X_\Delta(N)$ for such $N$ in Table
\ref{lowgenus}. Also the other cases can be found in Table
\ref{notrigonal}.

First applying Abramovich' bound we get the following result:
\begin{Lem}
The modular curves $X_{\Delta_i^\dag}$ and $X_{\Delta_i^\ddag}$ in
Table \ref{lowgenus},\ref{notrigonal} are not hyperelliptic.
\end{Lem}

\begin{Rk}
The notations ${\Delta_i}^\dag$ and ${\Delta_i}^\ddag$ in the
tables mean that Abramovich's bound doesn't hold for
$X_{{\Delta_i}^\dag}(N)$ and $X_{{\Delta_i}^\ddag}(N)$ when
$\operatorname{Gon}(X_{\Gamma_{{\Delta_i}^\dag}(N)})\leq 2$ and
$\operatorname{Gon}(X_{\Gamma_{{\Delta_i}^\ddag}(N)})\leq 3$
respectively.

\end{Rk}

Now we prove that $X_{\Delta_1}(21)$ is a hyperelliptic curve in
two different ways.

\noindent{\bf Proof 1.} The space $S_{\Delta_1}^2(21)$ is of
dimension $3$ and from \cite{St} we can get a basis consisting
three newforms. We write down the basis as follows:
\begin{eqnarray*}
f_1 &=& q - q^2 + q^3 - q^4 - 2q^5 - q^6 - q^7 + 3q^8 + q^9 +
2q^{10} + \cdots\\
f_2 &=& q - q^3 - 2q^4 + 2q^6 - 2q^7 + 4q^{10} + 2q^{11} + q^{13}
- 2q^{14} -\cdots\\
f_3 &=& 2q^2 - q^3 - 2q^4 - 2q^5 + q^7 + q^9 + 4q^{10} + 2q^{11} +
q^{13} - \cdots.
\end{eqnarray*}
By using the computer algebra system MAPLE we get a quadratic
generator of the ideal $I(\overline{X_{\Delta_1}(21)})$ as
follows:
$$Q\, :\, x_1^2-x_2^2-x_3^2+x_2x_3$$
where we obtain the relation $Q(f_1,f_2,f_3)=0$ by assigning $x_i$
to $f_i.$ But it means that $X_{\Delta_1}(21)$ is hyperelliptic by
the above criterion.\qed

\noindent{\bf Proof 2.} The authors \cite{J-K} proved that
$X_1(21)$ is bielliptic, i.e., it admits a map of degree $2$ to an
elliptic curve, and all the bielliptic involutions on $X_1(21)$
are $W_3=\begin{pmatrix} 9&-4\\21&-9\end{pmatrix}$ and $[8]W_3$
where $[a]$ denote the automorphism of $X_1(N)$ represented by
$\gamma\in\Gamma_0(N)$ such that $\gamma\equiv\sm a&*\\0&*\esm\md
N.$ Note that for bielliptic curves of genus $5$ all bielliptic
involutions commute with each other [Sch, Lemma 4.4]. Let $G$ be
the group generated by the two bielliptic involutions of
$X_1(21).$ Then we can determine the genus of the quotient
$G\backslash X_1(21)$ by the four-group rule \cite{F} as follows:
$$g(X_1(21))=g(W_3\backslash X_1(21))+g([8]W_3\backslash X_1(21))
+g([8]\backslash X_1(21))-2g(G\backslash X_1(21)).$$ Thus
$G\backslash X_1(21)$ is rational, and hence we get a Galois
covering $X_1(21)\to{\mathbb P}^1$ with Galois group $G.$ Since
$[8]\backslash X_1(21)$ is the same as $X_{\Delta_1}(21),$ we can
conclude that $X_{\Delta_1}(21)$ is hyperelliptic.\qed

To determine all the other curves $X_\Delta(N)$ to be not
hyperelliptic it suffices to consider $X_\Delta(N)$ for the
maximal subgroups $\Delta.$ For example, the modular curve
$X_{\Delta_1}(30)$ is of genus $5$ and has a basis of
$S_{\Delta_1}^2(30)$ which consists of two old forms and three new
forms as follows:
\begin{eqnarray*}
f_1 &=& q - q^2 - q^3 - q^4 + q^5 + q^6 + 3q^8 + q^9 - q^{10} -  \cdots\\
f_2 &=& q^2 - q^4 - q^6 - q^8 + q^{10} + q^{12} + 3q^{16} + q^{18} - q^{20} - \cdots\\
f_3 &=& q - q^2 + q^3 + q^4 - q^5 - q^6 - 4q^7 - q^8 + q^9 +
q^{10} +  \cdots\\
f_4 &=& q - q^4 - 2q^5 + q^6 - q^9 - q^{10} + 2q^{11} + 2q^{14} +
q^{15} + \cdots\\
f_5 &=& q^2 - q^3 + q^5 - 2q^7 - q^8 - 2q^{10} + q^{12} + 6q^{13}
+ 2q^{15} - \cdots.
\end{eqnarray*}
By using MAPLE we get three quadratic generators of
$I(\overline{X_{\Delta_1}(30)})$ as follows:
$$\left\{
\begin{array}{ll}
    x_4^2-x_5^2-x_1x_3+2x_2x_3-4x_4x_5,\\
    x_3^2-2x_5^2+2x_1x_2-x_1x_3+2x_2x_3-4x_4x_5,\\
    x_1^2+4x_2^2-2x_5^2+2x_1x_2-x_1x_3+2x_2x_3-4x_4x_5.\\
\end{array}%
\right.$$ This means that $X_{\Delta_1}(30)$ is not hyperelliptic.
Case by case we calculate the quadratic generators of
$I(\overline{X_\Delta(N)})$ for maximal subgroups $\Delta$ in
Table \ref{lowgenus},\ref{notrigonal} and hence we can finish the
proof of Theorem \ref{hyperelliptic}.

\section{Trigonal modular curves}
In this section we determine all trigonal modular curves
$X_\Delta(N).$ Combining Theorem \ref{hyperelliptic} with
Proposition \ref{d-gonal} it suffices to consider the modular
curves $X_\Delta(N)$ of $g_\Delta(N)\geq 5$ in the Table
\ref{lowgenus},\ref{trigonal2} which contain all the intermediate
modular curves between $X_1(N)$ and $X_0(N)$ such that $X_0(N)$ is
trigonal.

Applying Abramovich's bound we get the following result.
\begin{Lem}\label{abramovich}
All the modular curves $X_{\Delta_i^\ddag}(N)$ in Table
\ref{lowgenus},\ref{trigonal2} are not trigonal.
\end{Lem}

We make use of the method due to Hasegawa and Shimura \cite{H-S}.

\begin{Thm}{\rm (Petri's Theorem)} Let $X$ be a canonical curve of
genus $g\geq 4$ defined over an algebraically closed field. Then
the ideal $I(X)$ of $X$ is generated by some quadratic
polynomials, unless $X$ is trigonal or isomorphic to a smooth
plane quintic curve, in which cases it is generated by some
quadratic and (at least one) cubic polynomials.
\end{Thm}

Let $X_\Delta(N)$ be of genus $g_\Delta(N)\geq 5$ and
$\{f_1,f_2,\dots,f_g\}$ a basis of $S_\Delta^2(N).$ Then to obtain
a minimal generating system of the ideal
$I(\overline{X_\Delta(N)})$, we have only to compute the relations
of the $f_if_j$ and the $f_if_jf_k(1\leq i,j,k\leq g),$ and to
eliminate those cubic relations arising from quadratic relations.
By Petri's Theorem $X_\Delta(N)$ is trigonal if and only if it is
not isomorphic to a smooth plane quintic curve and a minimal
generating system of $I(\overline{X_\Delta(N)})$ contains a cubic
polynomial. Let $Q_1,\dots,Q_{(g-2)(g-3)/2}$ be a system of
quadratic generators of $I(\overline{X_\Delta(N)}).$ Since there
are $(g-3)(g^2+6g-10)/6$ linear relations among the $f_if_jf_k,$
the number of cubic generators among the minimal generating system
is given by
$$\frac{(g-3)(g^2+6g-10)}{6}-\dim L'$$
where $L'$ is the subspace generated by $x_iQ_j\,\,(1\leq i\leq
g;\,1\leq j\leq (g-2)(g-3)/2).$ Thus $X_\Delta(N)$ is trigonal
only if the above difference is non-zero.

\begin{Exm}
The curve $X_{\Delta_1}(32)$ is of genus $5$ and not
hyperelliptic. By the exact same method as in the computation of
$X_{\Delta_1}(30)$ (See \S 2) we get three quadratic generators of
$I(\overline{X_{\Delta_1}(32)})$ as follows:
$$\left\{
\begin{array}{ll}
x_1^2+x_2^2+x_3^2+8x_5^2+2x_2x_3+4x_2x_4-4x_2x_5-8x_4x_5\\
-x_2x_3-x_2x_4-x_2x_5-x_3x_4+x_3x_5\\
x_4^2-x_5^2+x_2x_5+x_3x_4+2x_4x_5.
\end{array}\right.$$ By a simple
calculation we find that the dimension of $L'$ is exactly $15,$ it
follows that there are no essential cubic generators. Therefore
$X_{\Delta_1}(30)$ is not trigonal.
\end{Exm}

Following the same method as in the above example we calculate the
remaining cases to get the following result.

\begin{Thm}\label{trigonal}
The modular curve $X_\Delta(N)$ is trigonal if and only if it is
of genus $g_\Delta(N)\leq 2$ or not hyperelliptic of
$g_\Delta(N)=3,4.$ Equivalently all the curves $X_\Delta(N)$ of
genus $g_\Delta(N)\leq 4$ in Table \ref{lowgenus} except
$X_{\Delta_1}(21).$
\end{Thm}

\begin{center}
{\bf Acknowledgment}
\end{center}
We thank to Andreas Schweizer for suggesting the second proof of
Theorem \ref{hyperelliptic}.

\section*{appendix}

\begin{center}

\begin{longtable}{|c|l|c|}
\caption{List of $X_\Delta(N)$ and their genera $g_\Delta(N)$ when
$X_0(N)$ are of genus $g_0(N)\leq 2$.}\label{lowgenus}\\
\hline $N$ & $\{\pm 1\} \subsetneq \Delta \subsetneq (\Z/N\Z)^*$ &
$g_\Delta(N)$
 \\ \hline
 $1\le N \le 12$ & $-$ & $-$
 \\ \hline
 $13$ & $\Delta_1=\{ \pm 1, \pm 5\}$ & 0
 \\ \hline
 $13$ & $\Delta_2=\{ \pm 1, \pm 3, \pm 4 \}$ & 0
 \\ \hline
 $14$ & $-$ & $-$
 \\ \hline
 $15$ & $\Delta_1=\{ \pm 1, \pm 4\}$ & 1
 \\ \hline
 $16$ & $\Delta_1=\{ \pm 1, \pm 7\}$ & 0
 \\ \hline
 $17$ & $\Delta_1=\{\pm 1,\pm 4\}$ & $1$
 \\ \hline
 $17$ & $\Delta_2=\{\pm 1,\pm 2,\pm 4,\pm 8\}$ & $1$
 \\ \hline
 $18$ & $-$ & $-$
 \\ \hline
 $19$ & $\Delta_1=\{\pm 1,\pm 7,\pm 8\}$ & $1$
 \\ \hline
 $20$ & $\Delta_1=\{ \pm 1, \pm 9\}$ & 1
 \\ \hline
 $21$ & $\Delta_1=\{\pm 1,\pm 8\}$ & $3$
 \\ \hline
 $21$ & $\Delta_2=\{\pm 1,\pm 4,\pm 5\}$ & $1$
 \\ \hline
 $22$ & $-$ & $-$
 \\ \hline
 $23$ & $-$ & $-$
 \\ \hline
 $24$ & $\Delta_1=\{\pm 1,\pm 5\}$ & $3$
 \\ \hline
 $24$ & $\Delta_2=\{\pm 1,\pm 7\}$ & $3$
 \\ \hline
 $24$ & $\Delta_3=\{\pm 1,\pm 11\}$ & $1$
 \\ \hline
 $25$ & $\Delta_1=\{\pm 1,\pm 7\}$ & $4$
 \\ \hline
 $25$ & $\Delta_2=\{\pm 1,\pm 4,\pm 6,\pm 9,\pm 11\}$ & $0$
 \\ \hline
 $26$ & $\Delta_1=\{\pm 1,\pm 5\}$ & $4$
 \\ \hline
 $26$ & $\Delta_2=\{\pm 1,\pm 3,\pm 9\}$ & $4$
 \\ \hline
 $27$ & $\Delta_1=\{\pm 1,\pm 8,\pm 10\}$ & $1$
 \\ \hline
 $28$ & $\Delta_1=\{\pm 1,\pm 13\}$ & $4$
 \\ \hline
 $28$ & $\Delta_2=\{\pm 1,\pm 3,\pm 9\}$ & $4$
 \\ \hline
 $29$ & $\Delta_1^\dag=\{\pm 1,\pm 12\}$ & $8$
 \\ \hline
 $29$ & $\Delta_2=\{\pm 1,\pm 4,\pm 5,\pm 6,\pm 7,\pm 9,\pm 13\}$ & $4$
 \\ \hline
 $31$ & $\Delta_1=\{\pm 1,\pm 5,\pm 6\}$ & $6$
 \\ \hline
 $31$ & $\Delta_2=\{\pm 1,\pm 2,\pm 4,\pm 8,\pm 15\}$ & $6$
 \\ \hline
 $32$ & $\Delta_1=\{\pm 1,\pm 15\}$ & $5$
 \\ \hline
 $32$ & $\Delta_2=\{\pm 1,\pm 7,\pm 9,\pm 15\}$ & $1$
 \\ \hline
 $36$ & $\Delta_1^\dag=\{\pm 1,\pm 17\}$ & $7$
 \\ \hline
 $36$ & $\Delta_2=\{\pm 1,\pm 11,\pm 13\}$ & $3$
 \\ \hline
 $37$ & $\Delta_1^\ddag=\{\pm 1,\pm 6\}$ & $16$
 \\ \hline
 $37$ & $\Delta_2^\dag=\{\pm 1,\pm 10,\pm 11\}$ & $10$
 \\ \hline
 $37$ & $\Delta_3=\{\pm 1,\pm 6,\pm 8,\pm 10,\pm 11,\pm 14\}$ & $4$
 \\ \hline
 $37$ & $\Delta_4=\{\pm 1,\pm 3,\pm 4,\pm 7,\pm 9,\pm 10,\pm 11,\pm 12,\pm 16\} $ & $4$
 \\ \hline
 $49$ & $\Delta_1^\ddag=\{\pm 1,\pm 18,\pm 19\}$ & $19$
 \\ \hline
 $49$ & $\Delta_2=\{\pm 1,\pm 6,\pm 8,\pm 13,\pm 15,\pm 20,\pm 22\}$ & $3$
 \\ \hline
 $50$ & $\Delta_1^\ddag=\{\pm 1,\pm 7\}$ & $22$
 \\ \hline
 $50$ & $\Delta_2=\{\pm 1,\pm 9,\pm 11,\pm 19,\pm 21\}$ & $4$
 \\ \hline
\end{longtable}
\end{center}

\begin{center}
\begin{longtable}{|c|l|c|}
\caption{List of $X_\Delta(N)$ and their genera $g_\Delta(N)$ when
$X_0(N)$ are hyperelliptic and $g_0(N)>2$.}\label{notrigonal}\\

\hline $N$ & $\{ \pm 1 \} \subsetneq \Delta \subsetneq (\Z/N\Z)^*$
& $g_\Delta(N)$
 \\ \hline
 $30$ & $\Delta_1=\{\pm 1,\pm 11\}$ & $5$
 \\ \hline
 $33$ & $\Delta_1^\dag=\{\pm 1,\pm 10\}$ & $11$
 \\ \hline
 $33$ & $\Delta_2=\{\pm 1,\pm 2,\pm 4,\pm 8,\pm 16\}$ & $5$
 \\ \hline
 $35$ & $\Delta_1^\dag=\{\pm 1,\pm 6\}$ & $13$
 \\ \hline
 $35$ & $\Delta_2=\{\pm 1,\pm 11,\pm 16\}$ & $9$
 \\ \hline
 $35$ & $\Delta_3=\{\pm 1,\pm 6, \pm 8, \pm 13\}$ & $7$
 \\ \hline
 $35$ & $\Delta_4=\{\pm 1,\pm 4,\pm 6,\pm 9,\pm 11,\pm 16\}$ & $5$
 \\ \hline
 $39$ & $\Delta_1^\dag=\{\pm 1,\pm 14\}$ & $17$
 \\ \hline
 $39$ & $\Delta_2^\dag=\{\pm 1,\pm 16,\pm 17\}$ & $9$
 \\ \hline
 $39$ & $\Delta_3=\{\pm 1,\pm 5,\pm 8, \pm 14\}$ & $9$
 \\ \hline
 $39$ & $\Delta_4=\{\pm 1,\pm 4,\pm 10,\pm 14,\pm 16,\pm 17\}$ & $5$
 \\ \hline
 $40$ & $\Delta_1^\dag=\{\pm 1,\pm 31\}$ & $9$
 \\ \hline
 $40$ & $\Delta_2^\dag=\{\pm 1,\pm 9\}$ & $13$
 \\ \hline
 $40$ & $\Delta_3^\dag=\{\pm 1,\pm 11\}$ & $13$
 \\ \hline
 $40$ & $\Delta_4=\{\pm 1,\pm 9,\pm 11,\pm 19\}$ & $5$
 \\ \hline
 $40$ & $\Delta_5=\{\pm 1,\pm 3,\pm 9,\pm 13\}$ & $7$
 \\ \hline
 $40$ & $\Delta_6=\{\pm 1,\pm 7,\pm 9,\pm 17\}$ & $7$
 \\ \hline
 $41$ & $\Delta_1^\dag=\{\pm 1,\pm 9\}$ & $21$
 \\ \hline
 $41$ & $\Delta_2^\dag=\{\pm 1,\pm 3,\pm 9,\pm 14\}$ & $11$
 \\ \hline
 $41$ & $\Delta_3=\{\pm 1,\pm 4,\pm 10,\pm 16,\pm 18\}$ & $11$
 \\ \hline
 $41$ & $\Delta_4=\{\pm 1,\pm 2,\pm 4,\pm 5,\pm 8,\pm 9,\pm 10,\pm 16,\pm 18,\pm 20\}$ & $5$
 \\ \hline
 $46$ & $-$ & $-$
 \\ \hline
 $47$ & $-$ & $-$
 \\ \hline
 $48$ & $\Delta_1^\dag=\{\pm 1,\pm 7\}$ & $19$
 \\ \hline
 $48$ & $\Delta_2^\dag=\{\pm 1,\pm 17\}$ & $19$
 \\ \hline
 $48$ & $\Delta_3^\dag=\{\pm 1,\pm 23\}$ & $19$
 \\ \hline
 $48$ & $\Delta_4=\{\pm 1,\pm 11,\pm 13,\pm 23\}$ & $5$
 \\ \hline
 $48$ & $\Delta_5=\{\pm 1,\pm 7,\pm 17,\pm 23\}$ & $7$
 \\ \hline
 $48$ & $\Delta_6=\{\pm 1,\pm 5,\pm 19,\pm 23\}$ & $7$
 \\ \hline
 $59$ & $-$ & $-$
 \\ \hline
 $71$ & $\Delta_1^\dag=\{\pm 1,\pm 5,\pm 14,\pm 17,\pm 25\}$ & $36$
 \\ \hline
 $71$ & $\Delta_2^\dag=\{\pm 1,\pm 20,\pm 23,\pm 26,\pm 30,\pm 32,\pm 34\}$ & $26$
 \\ \hline
\end{longtable}
\end{center}

\begin{center}
\begin{longtable}{|c|l|c|}
\caption{List of $X_\Delta(N)$ and their genera $g_\Delta(N)$ when
$X_0(N)$ are trigonal but not
sub-hyperelliptic.}\label{trigonal2}\\

\hline $N$ & $\{ \pm 1 \} \subsetneq \Delta \subsetneq (\Z/N\Z)^*$
& $g_\Delta(N)$
 \\ \hline
 $34$ & $\Delta_1=\{\pm 1,\pm 13\}$ & $9$
 \\ \hline
 $34$ & $\Delta_2=\{\pm 1,\pm 9,\pm 13,\pm 15\}$ & $5$
 \\ \hline
 $38$ & $\Delta_1=\{\pm 1,\pm 7,\pm 11\}$ & $10$
 \\ \hline
 $43$ & $\Delta_1^\ddag=\{\pm 1,\pm 6,\pm 7\}$ & $15$
 \\ \hline
 $43$ & $\Delta_2=\{\pm 1,\pm 2,\pm 4,\pm 8,\pm 11,\pm 16,\pm 21,\pm 22\}$ & $9$
 \\ \hline
 $44$ & $\Delta_1^\ddag=\{\pm 1,\pm 21\}$ & $16$
 \\ \hline
 $44$ & $\Delta_2=\{\pm 1,\pm 5,\pm 7,\pm 9,\pm 19\}$ & $8$
 \\ \hline
 $45$ & $\Delta_1^\ddag=\{\pm 1,\pm 19\}$ & $21$
 \\ \hline
 $45$ & $\Delta_2=\{\pm 1,\pm 14,\pm 16\}$ & $9$
 \\ \hline
 $45$ & $\Delta_3=\{\pm 1,\pm 8,\pm 17,\pm 19\}$ & $11$
 \\ \hline
 $45$ & $\Delta_4=\{\pm 1,\pm 4,\pm 11,\pm 14,\pm 16,\pm 19\}$ & $5$
 \\ \hline
 $53$ & $\Delta_1^\ddag=\{\pm 1,\pm 23\}$ & $40$
 \\ \hline
 $53$ & $\Delta_2=\{\pm 1,\pm 4,\pm 6,\pm 7,\pm 9,\pm 10,\pm 11,\pm 13,\pm 15,\pm 16,\pm 17,\pm 24,\pm 25\}$ & $8$
 \\ \hline
 $54$ & $\Delta_1^\ddag=\{\pm 1,\pm 17,\pm 19\}$ & $10$
 \\ \hline
 $61$ & $\Delta_1^\ddag=\{\pm 1,\pm 11\}$ & $56$
 \\ \hline
 $61$ & $\Delta_2^\ddag=\{\pm 1,\pm 13,\pm 14\}$ & $36$
 \\ \hline
 $61$ & $\Delta_3^\ddag=\{\pm 1,\pm 3,\pm 9,\pm 20,\pm 27\}$ & $26$
 \\ \hline
 $61$ & $\Delta_4^\ddag=\{\pm 1,\pm 11,\pm 13,\pm 14,\pm 21,\pm 29\}$ & $16$
 \\ \hline
 $61$ & $\Delta_5=\{\pm 1,\pm 3,\pm 8,\pm 9,\pm 11,\pm 20,\pm 23,\pm 24,\pm 27,\pm 28\}$ & $12$
 \\ \hline
 $64$ & $\Delta_1^\ddag=\{\pm 1,\pm 31\}$ & $37$
 \\ \hline
 $64$ & $\Delta_2^\ddag=\{\pm 1,\pm 15,\pm 17,\pm 31\}$ & $13$
 \\ \hline
 $64$ & $\Delta_3=\{\pm 1,\pm 7,\pm 9,\pm 15,\pm 17,\pm 23,\pm 25,\pm 31\}$ & $5$
 \\ \hline
 $81$ & $\Delta_1^\ddag=\{\pm 1,\pm 26,\pm 28\}$ & $46$
 \\ \hline
 $81$ & $\Delta_2^\ddag=\{\pm 1,\pm 8,\pm 10,\pm 17,\pm 19,\pm 26,\pm 28,\pm 35,\pm 37\}$ & $10$
 \\ \hline
\end{longtable}
\end{center}

\end{document}